\documentclass{amsart}

\usepackage{xypic}
 \input xy
\xyoption{all}
\usepackage{epsfig}
\usepackage{amsthm}
\usepackage{amssymb}
\usepackage{amsmath}
\usepackage{amscd}

\newcommand{\bg}{\begin{equation}}
\newcommand{\ed}{\end{equation}}
\newcommand{\bga}{\begin{eqnarray}}
\newcommand{\eda}{\end{eqnarray}}
\newcommand{\pf}{\textbf{Proof:\ }}

\def\cbdu{\par{\raggedleft$\Box$\par}}

\newtheorem {Theorem}  {Theorem}

\numberwithin{Theorem}{section}

\newtheorem {Lemma}[Theorem]  {Lemma}
\newtheorem {Proposition}[Theorem]{Proposition}
\theoremstyle{definition}

\theoremstyle{remark}
\newtheorem{Remark}[Theorem]{\bf Remark}

\newtheorem {Corollary}[Theorem]{\bf Corollary}

\expandafter\chardef\csname pre amssym.def
at\endcsname=\the\catcode`\@ \catcode`\@=11
\def\undefine#1{\let#1\undefined}
\def\newsymbol#1#2#3#4#5{\let\next@\relax
 \ifnum#2=\@ne\let\next@\msafam@\else
 \ifnum#2=\tw@\let\next@\msbfam@\fi\fi
 \mathchardef#1="#3\next@#4#5}
\def\mathhexbox@#1#2#3{\relax
 \ifmmode\mathpalette{}{\m@th\mathchar"#1#2#3}%
 \else\leavevmode\hbox{$\m@th\mathchar"#1#2#3$}\fi}
\def\hexnumber@#1{\ifcase#1 0\or 1\or 2\or 3\or 4\or 5\or 6\or 7\or 8\or
 9\or A\or B\or C\or D\or E\or F\fi}

\font\teneufm=eufm10 \font\seveneufm=eufm7 \font\fiveeufm=eufm5
\newfam\eufmfam
\textfont\eufmfam=\teneufm \scriptfont\eufmfam=\seveneufm
\scriptscriptfont\eufmfam=\fiveeufm

\catcode`\@=\csname pre amssym.def at\endcsname

\newcounter{remark}
\usepackage{color}

\newcommand{\e}{\epsilon}

\newcommand{\R}{\mathbf{R}}

\def  \R   {{\mathbb R}}

\def  \12  {{\frac{1}{2}}}

\def\build#1_#2^#3{\mathrel{\mathop{\kern 0pt#1}\limits_{#2}^{#3}}}

\begin{document}

\title{Asymptotic Behavior of Solutions to the Liquid Crystals System in $\mathbb{R}^3$}

\author[Mimi Dai]{ Mimi Dai}
\address{Department of Mathematics, UC Santa Cruz, Santa Cruz, CA 95064,USA}
\email{mdai@ucsc.edu}
\author [Jie Qing]{Jie Qing}
\address{Department of Mathematics, UC Santa Cruz, Santa Cruz, CA 95064,USA}
\email{qing@ucsc.edu}
\author[Maria Schonbek] {Maria Schonbek}
\address{Department of MAthematics, UC Santa Cruz, Santa Cruz, CA 95064, USA}
\email{schonbek@ucsc.edu}




\begin{abstract}
 In this paper we study the large time behavior of  solutions to a nematic liquid crystals system in the whole space $\R^3$. The fluid under consideration has  constant  density  and small  initial data.  \end{abstract}

\maketitle

\section{Introduction}
In this paper we consider the asymptotic behavior of solutions to the simplified model of nematic liquid crystals (LCD) with constant density:

\begin{equation}\begin{split}\label{LCD}
u_t+u\cdot\nabla u+\nabla p =\nu\triangle u -\nabla\cdot(\nabla d\otimes\nabla d),\\
d_t+u\cdot\nabla d =\triangle d-f(d),\\
\nabla\cdot u =0.
\end{split}
\end{equation}
The equations are considered in $\R^3\times(0, T)$. Here
$p: \R^3\times[0,T]\to\mathbb{R}$ is the fluid pressure,  $
u: \R^3\times[0,T]\to\mathbb{R}^3$ is the fluid velocity and
$d: \R^3\times[0,T]\to\mathbb{R}^3$ is the direction field representing the alignment of the molecules. The constant $\nu >0$ stands for the viscosity coefficient.  Without loss of generality, by scaling, we can set $\nu =1$.  The force term $\nabla d\otimes\nabla d$ in the equation of the conservation of momentum denotes the $3\times 3$ matrix whose $ij$-th entry is given by $``\nabla_i d\cdot\nabla_j d"$ for $1\leq i,j\leq 3$. This force $\nabla d\otimes\nabla d$ is the stress tensor of the energy about the direction field $d$, where the   energy is given by:
$$
\frac 12 \int_{\R^3} |\nabla d|^2 dx + \int_{\R^3}F(d)dx
$$
where
$$
F(d)=\frac{1}{4\eta^2}(|d|^2-1)^2, \quad f(d) = \nabla F(d) = \frac{1}{\eta^2}(|d|^2-1)d,
$$
for a constant $\eta$ in this paper. We note that  $F(d)$ is the penalty term of the Ginzburg-Landau approximation of the original free energy of the direction field with unit length.

In this paper we consider the following  initial conditions:
\bg\label{initu}
u(x,0) =u_0(x), \ \ \ \nabla\cdot u_0 =0,
\ed
\bg\label{initd}
d(x,0) =d_0(x), \ \ \ |d_0(x)|= 1,
\ed
and
\bg\label{bd}
u_0\in H^1(\R^3), \ \ d_0-w_0\in H^2(\R^3),
\ed
with a fixed vector $w_0\in S^2$, i.e., $|w_0|=1$.

The flow of nematic liquid crystals can be treated as slow moving particles where the fluid velocity and the alignment of the particles influence each other. The hydrodynamic theory of liquid crystals  was established by Ericksen \cite{Er0, Er1} and Leslie \cite{Le0, Le1} in the 1960's.
As Leslie points out in his 1968 paper: ``liquid crystals are states of matter which are capable of flow, and in which the molecular arrangements  give rise to a preferred direction".  There is a vast literature on the hydrodynamic of  liquid crystal systems. For background we list  a few, with no intention to be complete:  \cite{EK, HKL, Kin, LL, LL2, LL1, Cal, CC, CC1, CDLL, Wu, Liu, JT}.  In particular, the asymptotic behavior of solutions to the flow of nematic liquid crystals was studied for bounded domains in \cite{LL, Wu}.  It was shown in \cite{Wu} that, with suitable initial conditions, the velocity converges to zero and the direction field converges to the steady solution
to the following equation
\begin{equation}\label{eq:steady}
\begin{cases}
-\Delta d+f(d)=0, x\in\Omega\\
d(x)=d_0(x), x\in\partial\Omega.
\end{cases}
\end{equation}

In  \cite{Wu}, Lemma 2.1 the \L ojasiewicz-Simon inequality  is used to derive the convergence when $\Omega$ is a bounded domain.
Lack of compactness considerations do not allow us to use similar arguments in the whole space $\R^3$.


In this paper we consider the asymptotic  behavior of the solutions  to \eqref{LCD} in the whole space $\R^3$.
subject to the additional condition on  the direction field which insures that the direction tends to a constant unit vector $w_0$,  as the space variable tends to infinity:
\begin{equation}\label{d0w0}
\lim_{|x|\to\infty} d_0(x) = w_0.
\end{equation}
This simplifies the situation and allows us to obtain the stability without needing
the Liapunov reduction and \L ojasiewicz-Simon inequality, since $w_0$ is a non-degenerate steady solution to \eqref{eq:steady}.

We start from the basic energy estimates (\ref{basic}) and Ladyzhenskaya estimates (\ref{lady}) \cite{LS,DQS} (see the extension to the whole space in appendix of this paper) for the system (\ref{LCD}). We then establish the convergence of the direction field $d$ to the constant steady solution $w_0$ based on Gagliardo-Nirenberg interpolation techniques. More precisely,  the convergence obtained is in $L^p(\R^3)$ for any $p>1$, with an algebraic decay rate of  $(1+t)^{-\frac{3}{2}(1-\frac{1}{p})}$. We then focus on the conservation of momentum equation in (\ref{LCD}).  We  apply the Fourier splitting technique  \cite{Sch, Sch1, Sch2} to obtain the decay of  the velocity $u$ with an algebraic decay rate of $(1+t)^{-\frac{1}{4}}$ in $L^2(\R^3)$ norm. This rate coincides with the decay rate of solutions to Navier-Stokes equations with a force decaying at a rate of $(1+t)^{-\frac{3}{4}}$ \cite{Sch}.

The existence of global regular solutions of (\ref{LCD}) with the initial and bounadry conditions has been established in \cite{LL} (in \cite{DQS} for nonconstant density) provided that the viscosity constant is large enough or initial data are small enough.  Based on the arguments in  \cite{DQS}  the existence of global regular solutions  of (\ref{LCD}), for small initial data,  is established in the appendix as follows:

\begin{Theorem}\label{exist} Let $u_0$ and $d_0$ satisfy (\ref{initu})-(\ref{bd}). Assume that $u_0\in H^1(\R^3)$ and $d_0- w_0 \in H^{2}(\R^3)\cap L^1(\R^3)$ for a unit vector $w_0$. There is a positive small number $\e_0$ such that if
\begin{equation}\label{smalldata}
\|u_0\|_{H^1(\R^3)}^2+\|d_0 - w_0\|_{H^2(\R^3)}^2\leq\epsilon_0,
\end{equation}
then the system (\ref{LCD}) has a classical solution $(u,p,d)$ in the time period $(0, T)$, for all $T>0$. That is, for some $\alpha\in(0,1)$
\begin{equation}\label{reg}\begin{split}
u\in C^{1+\alpha/2,2+\alpha}((0,T)\times\R^3)\\
\nabla p\in C^{\alpha/2,\alpha}((0,T)\times\R^3)\\
d\in C^{1+\alpha/2,2+\alpha}((0,T)\times\R^3).
\end{split}\end{equation}
And the solution $(u,p,d)$ satisfies the following basic energy estimate and higher order energy estimate (also called Ladyzhenskaya energy estimate in \cite{DQS} and \cite{LL})
\begin{align}\label{ex:basic}
&\int_{\R^3}|u|^2+|\nabla d|^2+2F(d)dx+2\int_0^T\int_{\R^3}|\nabla u|^2+|\Delta d-f(d)|^2dxdt\\
&\leq \|u_0\|_{L^2(\R^3)}^2+\|\nabla d_0\|_{L^2(\R^3)}^2\notag
\end{align}
\begin{align}\label{ex:lady}
&\int_{\R^3}|\nabla u|^2+|\Delta d|^2dx+\int_0^T\int_{\R^3}|\Delta u|^2+|\nabla\Delta d|^2dxdt\\
&\leq C(\|u_0\|_{H^1(\R^3)}^2+\|d_0 - w_0\|_{H^2(\R^3)}^2).\notag
\end{align}
Furthermore, the solution $d$ satisfies
\bg\label{ineq:d-w0L1time}
\int_{\R^3}|d(t)-w_0|dx\leq (C_0t+\int_{\R^3}|d_0-w_0|dx)e^{Ct}
\ed
with the constants $C_0$ and $C$ depending only on initial data and on $\eta$, respectively.
\end{Theorem}

For the smooth solution obtained in Theorem \ref{exist},  our main asymptotic result  is:

\begin{Theorem}\label{Mthm}
Let $(u,p,d)$ be smooth solution obtained in Theorem \ref{exist}. Assume additionally $u_0\in L^1(\R^3)$ and $d_0-w_0\in L^p(\R^3)$, for any $p \geq 1$ and a unit vector $w_0$. There exists a small number $\e_0>0$ such that if
\begin{equation}\label{smalldataH2}
\|u_0\|_{H^1(\R^3)}^2+\|d_0-w_0\|_{H^2(\R^3)}^2\leq\epsilon_0,
\end{equation}
then
\bg\label{ddecay}
\|d(\cdot,t)-w_0\|_{L^p(\R^3)}\leq C\|d_0-w_0\|_{L^p(\R^3)}(1+t)^{-\frac{3}{2}(1-\frac{1}{p})},
\ed
\bg\label{dh1decay}
\|\nabla (d(\cdot,t)-w_0)\|_{L^2(\R^3)}^2\leq C(1+t)^{-\frac{3}{4}},
\ed
\bg\label{udecay}
\|u(\cdot,t)\|_{L^2(\R^3)}^2\leq C(1+t)^{-\frac{1}{2}},
\ed
where the various constants $C$ only depend on initial data.
\end{Theorem}

The paper is organized as follows: in Section 2 we establish the decay for  the difference $d-w_0$, using the basic energy estimate (\ref{ex:basic}) and the Ladyzhenskaya energy estimate (\ref{ex:lady}). Combining the decay of $d-w_0$  and  Fourier splitting technique \cite{Sch}, in Section 3   we  obtain an algebraic decay for the velocity $u$ in $L^2(\R^3)$. In the appendix, we sketch a proof for the existence theorem \ref{exist}.

\bigskip

\section{Convergence of the direction field}

In this section  we  study the  $L^p$ decay of the direction field $d-w_0$ and  the decay for the first derivative. The first step is to  derive a uniform estimate of $d-w_0$ in $L^2(\R^3)$. This yields a uniform estimate for $d-w_0$ in $L^p(\R^3)$ for any $p\geq 1$. This $L^p$ estimate, is the basis  to establish the decay results.
For ease of reading we  state  the basic energy estimate and Ladyzhenskaya energy estimate satisfied by the smooth solution $(u,p,d)$ (see appendix for details),
\begin{align}\label{basic}
&\|u\|_{L^2}^2 + \|\nabla d\|_{L^2}^2 +\|F(d)\|_{L^1}+ 2\int_0^T\|\nabla u\|_{L^2}^2 + \|\Delta d-f(d)\|_{L^2}^2dt\\
&\leq \|u_0\|_{L^2}^2 + \|\nabla d_0\|_{L^2}^2\notag
\end{align}
\begin{equation}\label{lady}
\|\nabla u\|_{L^2}^2 + \|\Delta d\|_{L^2}^2 + \int_0^T\|\Delta u\|_{L^2}^2 + \|\nabla\Delta d\|_{L^2}^2dt \leq C(\|u_0\|_{H^1}^2 + \|d_0-w_0\|_{H^2}^2).
\end{equation}

In the sequel  we need  to use a Gagliardo-Nirenberg interpolation inequality. For completeness we recall from \cite{F} the inequality here
\begin{Proposition}\label{PPGN}\cite{F}
 Let $w\in W^{m,p}(\R^n)\cap L^q(\R^n )$, for $1\leq p\leq \infty$ and $1\leq q\leq \infty$. Then
\bg\label{GN}
\|D^kw\|_{L^r(\R^n )}\leq C\|D^mw\|_{L^p(\R^n)}^a\|w\|_{L^q(\R^n )}^{1-a}
\ed
for any integer $k\in[0,m-1]$, where
\bg\label{parameter}
\frac{1}{r}=\frac{k}{n}+a(\frac{1}{p}-\frac{m}{n})+(1-a)\frac{1}{q}
\ed
with $a\in[\frac{k}{m},1]$, either if $p=1$ or $p>1$ and $m-k-\frac{n}{p}\notin\mathcal{N}\cup\left\{0\right\}$, while $a\in[\frac{k}{m},1)$, if $p>1$ and $m-k-\frac{n}{p}\in\mathcal{N}\cup\left\{0\right\}$.
\end{Proposition}

\bigskip

\subsection{Uniform estimate of $d-w_0$ in $L^2(\R^3)$}\label{subsection3}

In this part, we show  that  the integrals\\
\[  \int_{\R^3}|d(x,t)-w_0|^2dx \; \,\mbox{and}\;\; \int_0^T\int_{\R^3}|\nabla(d(x,t)-w_0)|^2dxdt\]
are uniformly bounded by the initial data, applying  estimates for the solutions obtained in Theorem \ref{exist}. We have the following lemma,
\begin{Lemma}\label{Le:dl2}
Let $d$ be the solution obtained in Theorem \ref{exist}. There exists $\epsilon_0$   sufficiently small so that if  $\|u_0\|_{L^2(\R^3)}+\|\nabla d_0\|_{L^2(\R^3)} \leq \epsilon_0 $. Then
\bg\label{d-w0l2}
\int_{\R^3}|d(x,t)-w_0|^2dx+\int_0^T\int_{\R^3}|\nabla(d(x,t)-w_0)|^2dxdt\leq C,
\ed
where the constant $C$ depends on the initial data and the norm $\|d_0-w_0\|_{L^2(\R^3)}$.
\end{Lemma}
\pf
Since $|w_0|=1$ and $f(w_0)=\frac{1}{\eta^2}(|w_0|^2-1)w_0=0$, the second equation in (\ref{LCD}) can be expressed as
\begin{equation}\label{eq:rewritten}
(d-w_0)_t+u\cdot\nabla(d-w_0)=\Delta(d-w_0)-f(d)+f(w_0).
\end{equation}
Applying the mean value theorem for vector valued functions, we have
\bg\label{mean}
f(d)-f(w_0)=\left (\int_0^1Df(w_0+s(d-w_0))ds\right )\cdot (d-w_0),
\ed
where $Df$ denotes the Jacobian matrix of $f$.
Multiplying (\ref{eq:rewritten}) by $d-w_0$ yields
\begin{align}\label{eq:L2d-w0}
&\frac{1}{2}\frac{d}{dt}\int_{\R^3}|d-w_0|^2dx\\
&=-\int_{\R^3}[u\cdot\nabla(d-w_0)](d-w_0)dx\notag\\
&+\int_{\R^3}\Delta(d-w_0)\cdot(d-w_0)dx\notag\\
&-\int_{\R^3}(d-w_0)^T\left (\int_0^1Df(w_0+s(d-w_0))ds\right )(d-w_0)dx,\notag\\
&\equiv -I_4+I_5+I_6\notag.
\end{align}
The three terms $I_4$, $I_5$ and $I_6$ are estimated as follows:\\
\begin{align}\label{est:I4}
|I_4|&=|\lim_{R\to\infty}\int_{B_R}[u\cdot\nabla(d-w_0)](d-w_0)dx|\\
&=|\lim_{R\to\infty}\frac{1}{2}\int_{\partial B_R}|d-w_0|^2u\cdot nd\sigma|\notag\\
&\leq C\lim_{R\to\infty}\left(\int_{\partial B_R}|d-w_0|^4d\sigma\right)^{1/2}\left(\int_{\partial B_R}|u|^2d\sigma\right)^{1/2}\notag\\
&\leq C\lim_{R\to\infty}\left(\int_{\partial B_R}|d-w_0|d\sigma\right)^{1/2}\left(\int_{\partial B_R}|u|^2d\sigma\right)^{1/2}\notag.
\end{align}
Denote the set $A=\left\{R:\int_{\partial B_R}|d-w_0|d\sigma\geq M\right\}$, for a certain constant $M>0$. $A$ is closed and its complement $A^c$ is open.  Write
$$
\int_{\R^3}|d(t)-w_0|dx=\int_0^\infty\int_{\partial B_R}|d(t)-w_0|d\sigma dR.
$$
Recall that by Theorem (\ref{exist})  the last expression is bounded by the initial data, for any fixed time $t>0$. Thus by Chebyshev's inequality  we have
\bg\label{measureA}
\mu \left\{A\right\}\leq \frac{C}{M}
\ed
where $\mu$ denotes the measure of a set, and $C$ denotes a constant depending only on the initial data.
Since the energy estimate (\ref{basic}) implies that
$$
\int_{\R^3}|u|^2dx=\int_0^\infty\int_{\partial B_R}|u|^2d\sigma dR
$$
is bounded by initial data, there exists a sequence $\left\{R_i\right\}_{i=1}^\infty\subset A^c$ with $R_i\to\infty$ such that
\bg\label{limit:uL2}
\int_{\partial B_{R_i}}|u|^2d\sigma\to 0
\ed
Combining the inequalities (\ref{est:I4}), (\ref{measureA}) and (\ref{limit:uL2}), yields  for all $t>0$,
\bg\label{limit:I4}
I_4=0.
\ed

For  $I_5$, we have
\begin{equation}\label{est:I5}
I_5=\lim_{R\to\infty}[\int_{\partial B_R}\frac{\partial(d-w_0)}{\partial n}\cdot(d-w_0)d\sigma-\int_{B_R}|\nabla(d-w_0)|^2dx]
\end{equation}
The boundary term is estimated as  follows:
\begin{align}\label{I5boundary}
&\int_{\partial B_R}\frac{\partial(d-w_0)}{\partial n}\cdot(d-w_0)d\sigma\\
&\leq\left(\int_{\partial B_R}|\nabla(d-w_0)|^2d\sigma\right)^{1/2}
\left(\int_{\partial B_R}|d-w_0|^2d\sigma\right)^{1/2}\notag\\
&\leq C\left(\int_{\partial B_R}|\nabla(d-w_0)|^2d\sigma\right)^{1/2}
\left(\int_{\partial B_R}|d-w_0|d\sigma\right)^{1/2}\notag
\end{align}
From Theorem \ref{exist}, we have that $\|\nabla(d-w_0)\|_{L^2(\R^3)}$ is uniformly bounded for any fixed $t>0$, $\|d(t)-w_0\|_{L^1(\R^3)}$ has a time dependent bound. For the inequality (\ref{I5boundary}), we apply a similar argument used to derive  (\ref{limit:I4}) and obtain the existence of a  sequence $R_i$ approaching infinity satisfying     \bg\notag
\lim_{R_i\to\infty} \int_{\partial B_{R_i}}\frac{\partial(d-w_0)}{\partial n}\cdot (d-w_0)d\sigma\to 0
\ed
It follows then that
\begin{equation}\label{limit:I5}
I_5=-\int_{\R^3}|\nabla(d-w_0)|^2dx.
\end{equation}

With respect to $I_6$, let $\tilde d=w_0+s(d-w_0)$ and $z=d-w_0$. Using the definition of $f(d)=\frac{1}{\eta^2}(|d|^2-1)d$, a straightforward calculation yields
\begin{align}\notag
&(d-w_0)^T\left (\int_0^1Df(w_0+s(d-w_0))ds\right )(d-w_0)=z^T\int_0^1Df(\tilde d)dsz\notag\\
&=\int_0^1[(2\tilde d_1^2+|\tilde d|^2-1)z_1^2+(2\tilde d_2^2+|\tilde d|^2-1)z_2^2+(2\tilde d_3^2+|\tilde d|^2-1)z_3^2\notag\\
&+4\tilde d_1\tilde d_2z_1z_2+4\tilde d_1\tilde d_3z_1z_3+4\tilde d_2\tilde d_3z_2z_3]ds\notag\\
&=\int_0^1[2\tilde d_1^2z_1^2+2\tilde d_2^2z_2^2+2\tilde d_3^2z_3^2+
4\tilde d_1\tilde d_2z_1z_2+4\tilde d_1\tilde d_3z_1z_3\notag\\
&+4\tilde d_2\tilde d_3z_2z_3]+(|\tilde d|^2-1)(z_1^2+z_2^2+z_3^2)ds\notag\\
&=\int_0^12[\tilde d_1z_1+\tilde d_2z_2+\tilde d_3z_3]^2+(|\tilde d|^2-1)|z|^2ds\notag\\
&=\int_0^12[\tilde d\cdot z]^2+(|\tilde d|^2-1)|z|^2ds\notag.
\end{align}
In the above equation, the third and forth equality comes from regrouping  terms and completing a perfect square.
Thus, $I_6$ can be written as
\begin{align}\label{est:I6}
I_6&=-\frac{1}{\eta^2}\int_{\R^3}\int_0^12[(w_0+s(d-w_0))\cdot(d-w_0)]^2\\
&+(|w_0+s(d-w_0)|^2-1)|d-w_0|^2dsdx\notag\\
&=-\frac{1}{\eta^2}\int_{\R^3}2[w_0\cdot(d-w_0)]^2+3w_0\cdot(d-w_0)|d-w_0|^2+|d-w_0|^4dx\notag\\
&=-\frac{1}{\eta^2}\int_{\R^3}2[w_0\cdot(d-w_0)+\frac{3}{4}|d-w_0|^2]^2-\frac{1}{8}|d-w_0|^4dx\notag\\
&\leq \frac{1}{8\eta^2}\int_{\R^3}|d-w_0|^4dx\notag.
\end{align}

Combining   (\ref{eq:L2d-w0}) and the inequalities (\ref{limit:I4}), (\ref{limit:I5}), and (\ref{est:I6}) gives
\bg\label{ml2:d-w0}
\frac{1}{2}\frac{d}{dt}\int_{\R^3}|d-w_0|^2dx+\int_{\R^3}|\nabla(d-w_0)|^2dx\leq \frac{1}{8\eta^2}\int_{\R^3}|d-w_0|^4dx.
\ed\\

The right hand side of the inequality (\ref{ml2:d-w0}) can be estimated as
\begin{align}\notag
\int_{\R^3}|d-w_0|^4dx&\leq (\int_{\R^3}|d-w_0|^2dx)^{1/2}(\int_{\R^3}|d-w_0|^6dx)^{1/2}\\
&\leq C(\int_{\R^3}|d-w_0|^2dx)^{1/2}(\int_{\R^3}|\nabla(d-w_0)|^2dx)^{3/2}\notag\\
&\leq C\int_{\R^3}|d-w_0|^2dx\left(\int_{\R^3}|\nabla(d-w_0)|^2dx\right )^{2}+\frac{1}{2}\int_{\R^3}|\nabla(d-w_0)|^2dx\notag.
\end{align}
 Gagliardo-Nirenberg interpolation inequality (Proposition \ref{PPGN}) yields
\bg\notag
\|\nabla(d-w_0)\|_{L^2(\R^3)}^2\leq C\|d-w_0\|_{L^2(\R^3)}\|\Delta(d-w_0)\|_{L^2(\R^3)}.
\ed
Combining the last two inequalities with (\ref{ml2:d-w0}) gives
\begin{align}\label{mml2:d-w0}
&\frac{d}{dt}\int_{\R^3}|d-w_0|^2dx+\int_{\R^3}|\nabla(d-w_0)|^2dx\\
&\leq C\left(\int_{\R^3}|d-w_0|^2dx\right)^2\int_{\R^3}|\Delta(d-w_0)|^2dx\notag.
\end{align}
Denote $\phi(t)=\int_{\R^3}|d(t)-w_0|^2dx$. Then
\bg\notag
\frac{d\phi}{\phi^2}\leq C\int_{\R^3}|\Delta(d-w_0)|^2dxdt.
\ed
Integrating the last inequality over $[0,t]$ yields
\bg\notag
-\frac{1}{\phi(t)}+\frac{1}{\phi(0)}\leq C\int_0^t\int_{\R^3}|\Delta(d-w_0)|^2dxdt.
\ed
Thus,
\bg\notag
\phi(t)\leq\frac{\phi(0)}{1-C\phi(0)\int_0^t\int_{\R^3}|\Delta(d-w_0)|^2dxdt}.
\ed
From the basic energy estimate (\ref{basic}),  and the hypothesis we have
$$\int_0^t\int_{\R^3}|\Delta(d-w_0)|^2dxdt\leq  \|u_0\|_{L^2}^2+\|\nabla d_0\|_{L^2}^2\leq \epsilon_0.$$ Assume that  $\epsilon$ is so small that $C\epsilon_0 \phi(0) < 1/2$, then
$$C\phi(0)\int_0^t\int_{\R^3}|\Delta(d-w_0)|^2dxdt<\frac{1}{2}.
$$
Hence  for any $t>0$,
$$
\phi(t)\leq 2\phi(0),
$$
that is,
\bg\label{globall2}
\int_{\R^3}|d(t)-w_0|^2dx\leq 2\int_{\R^3}|d_0-w_0|^2dx.
\ed\\

Due to the estimates (\ref{globall2}) and (\ref{mml2:d-w0}), we have
\bg\notag
\frac{d}{dt}\int_{\R^3}|d-w_0|^2dx+\int_{\R^3}|\nabla(d-w_0)|^2dx
\leq C\int_{\R^3}|\Delta(d-w_0)|^2dx,
\ed
where the constant $C$ only depends on the initial data. Integrating over $[0,t]$,  by the basic energy inequality  (\ref{basic}) it follows that
\begin{align}\notag
&\int_{\R^3}|d(t)-w_0|^2dx+\int_0^t\int_{\R^3}|\nabla(d-w_0)|^2dxdt\\
&\leq \int_{\R^3}|d_0-w_0|^2dx+C\int_0^t\int_{\R^3}|\Delta(d-w_0)|^2dxdt\notag\\
&\leq C\notag.
\end{align}
Thus,
\bg\label{l2spacetime}
\int_0^t\int_{\R^3}|\nabla(d-w_0)|^2dxdt\leq C,
\ed
where the constant $C$ only depends on initial data. This completes the proof of the lemma.
\cbdu

\bigskip

The following   auxiliary estimate  shows  that,   provided the  initial data is small enough, the norm $\|d(\cdot,t)-w_0\|_{L^\infty(\R^3)}$ will be  as small as necessary.  This smallness  yields  that   $|d|$ will be close to $1$, for all time.
\begin{Lemma}\label{Le:linfty}
Let $d$ be the solution obtained in Theorem \ref{exist}.  Then
\bg\label{dw0}
\|d(\cdot,t)-w_0\|_{L^\infty(\R^3)}\leq C\|\nabla d(t)\|_{L^2(\R^3)}^{1/2}\|\Delta d(t)\|_{L^2(\R^3)}^{1/2},
\ed
where $C$ is an absolute constant.
\end{Lemma}
\pf
Let $\mathcal F$ denote the Fourier transform. By Lemma \ref{Le:dl2}, we can take  the Fourier transform of $d-w_0$
\begin{align}\label{primo}
&\|d(\cdot,t)-w_0\|_{L^\infty(\R^3)} \leq \int_{\R^3}|\mathcal F(d-w_0)|d\xi\\
&=\int_{|\xi|\leq\lambda}\mathcal F(d-w_0)|d\xi+\int_{|\xi|\geq\lambda}\mathcal| F(d-w_0)|d\xi\notag\\
&=\int_{|\xi|\leq\lambda}\frac{1}{|\xi|}\cdot|\xi|\mathcal |F(d-w_0)|d\xi+\int_{|\xi|\geq\lambda}\frac{1}{|\xi|^2}\cdot|\xi^2||\mathcal F(d-w_0)|d\xi\notag\\
&\leq (\int_{|\xi|\leq\lambda}\frac{1}{|\xi|^2}d\xi)^{\frac{1}{2}}(\int_{|\xi|\leq\lambda}|\xi||\mathcal F(d-w_0)|^2d\xi)^{\frac{1}{2}}+(\int_{|\xi|\geq\lambda}\frac{1}{|\xi|^4}d\xi)^{\frac{1}{2}}(\int_{|\xi|\geq\lambda}|\xi^2\mathcal F(d-w_0)|^2d\xi)^{\frac{1}{2}}\notag\\
&\leq C(\int_{r\leq\lambda}\frac{r^2}{r^2}dr)^{\frac{1}{2}}(\int_{\R^3}|\nabla(d-w_0)|^2dx)^{\frac{1}{2}}+
C(\int_{r\geq\lambda}\frac{r^2}{r^4}dr)^{\frac{1}{2}}(\int_{\R^3}|\Delta(d-w_0)|^2dx)^{\frac{1}{2}}\notag\\
&\leq C\lambda^{\frac{1}{2}}\|\nabla d(t)\|_{L^2(\R^3)}+C\lambda^{-\frac{1}{2}}\|\Delta d(t)\|_{L^2(\R^3)}\notag.
\end{align}
 To find a $\lambda$ that minimizes the right hand side of the last inequality,
 take the derivative in $\lambda$ and set the right hand side equal to zero:

$$\lambda^{\frac{1}{2}}\|\nabla d(t)\|_{L^2(\R^3)}=\lambda^{-\frac{1}{2}}\|\Delta d(t)\|_{L^2(\R^3)},$$
yielding
$$
\lambda=\|\Delta d(t)\|_{L^2(\R^3)}/\|\nabla d(t)\|_{L^2(\R^3)}.
$$
Thus using this $\lambda $ in (\ref{primo}) gives inequality (\ref{dw0}), and the proof of the Lemma is complete.

\cbdu

\begin{Corollary}
Suppose the initial data $\|u_0\|_{H^1} + \|d_0-w_0\|_{H^2}$ are small enough. Then, $|d(x,t)|\geq \frac{1}{2}$.
\end{Corollary}
\pf It follows combining (\ref{basic}), (\ref{lady}), (\ref{dw0}) since $|w_0|=1$. \cbdu

\bigskip

\subsection{Uniform estimate of $d-w_0$ in $L^p$ with any $p\geq 1$}

Here we show provided the data is small enough,  all the $L^p$ norms of $d-w_0$  are bounded.
\begin{Lemma}\label{Le:dlp}
Let $d$ be the solution obtained in Theorem \ref{exist}. There exist  $\lambda_p $, depending on $p$,  so that if  $\|u_0\|_{H^1(\R^3)}^2+\|d_0-w_0\|_{H^2(\R^3)}^2\leq \lambda_p$, then  for  $p> 1$
\begin{align}\label{d-w0lp}
&\frac{1}{p}\int_{\R^3}|d(x,t)-w_0|^pdx+\frac{2(p-1)}{p^2}\int_0^T\int_{\R^3}|\nabla|d(x,t)-w_0|^{p/2}|^2dxdt\\
&\leq C_p\int_{\R^3}|d_0-w_0|^pdx\notag,
\end{align}
where the constant $C_p$ depends on $p$ and  $\lambda_p$. And for $p=1$ we have
\begin{align}\notag
\int_{\R^3}|d-w_0|dx&\leq\int_{\R^3}|d_0-w_0|dx+\int_0^T\int_{\R^3}|\nabla(d-w_0)|^2+|\Delta(d-w_0)|^2dxdt\notag\\
&\leq C\notag.
\end{align}

\end{Lemma}
\pf
Recall that since $f(w_0)=0$, we have that the direction equation can be rewritten as
 \bg\notag
(d-w_0)_t+u\cdot\nabla(d-w_0)=\Delta(d-w_0)-f(d)+f(w_0).
\ed
Multiplying the last equation by  $(d-w_0)|d-w_0|^{p-2}$, for any $p \geq 2$ (or alternatively  by  $(d- w_0) /(|d - w_0|+ \epsilon)^{2-p}$, when $p\in [1, 2)$, and  letting $\epsilon\to 0$) yields
\begin{align}\label{eq:Lpd-w0}
&\frac{1}{p}\frac{d}{dt}\int_{\R^3}|d-w_0|^pdx\\
&=-\int_{\R^3}[u\cdot\nabla(d-w_0)](d-w_0)|d-w_0|^{p-2}dx\notag\\
&+\int_{\R^3}\Delta(d-w_0)(d-w_0)|d-w_0|^{p-2}dx\notag\\
&-\int_{\R^3}2[w_0\cdot(d-w_0)+\frac{3}{4}|d-w_0|^2]^2|d-w_0|^{p-2}-\frac{1}{8}|d-w_0|^{p+2}dx\notag\\
&\equiv I_7+I_8+I_9\notag,
\end{align}
where the $I_9$ was obtained similarly as in the previous calculation for $I_6$.  We  estimate $I_7$, $I_8$ and $I_9$  as follows:\\
Integrating by parts over ball $B_R$ gives
\begin{equation}\notag
I_7=\lim_{R\to\infty}[\int_{\partial B_R}|d-w_0|^pu\cdot nd\sigma-(p-1)\int_{B_R}[u\cdot\nabla(d-w_0)](d-w_0)|d-w_0|^{p-2}dx].
\end{equation}
It implies that
\begin{align}\label{eq:I7}
pI_7&=\lim_{R\to\infty}\int_{\partial B_R}|d-w_0|^pu\cdot nd\sigma\\
&\leq \lim_{R\to\infty}\left(\int_{\partial B_R}|d-w_0|^2d\sigma\right)^{1/2}\left(\int_{\partial B_R}|u|^2d\sigma\right)^{1/2}\notag
\end{align}
for any $p\geq 1$, where we used  that $|d-w_0|\leq C$. \\
By Lemma \ref{Le:dl2} we know that  $\int_{\R^3}|d-w_0|^2dx$ is bounded, and $\int_{\R^3}|u|^2dx$ is bounded from the energy estimate (\ref{basic}). Thus, for the inequality (\ref{eq:I7}),  using  arguments  similar to the ones  applied to derive  the convergence (\ref{limit:I4}),  will yield
\bg\notag
\int_{\partial B_{R_i}}|d-w_0|^pu\cdot nd\sigma\to 0
\ed
 for an appropriate  sequence $R_i\to\infty$,  for any $p\geq 1$. Thus,
\bg\label{limit:I7}
I_7 =0.
\ed
For $I_8$, integrating by parts over $B_R$ yields
\begin{align}\label{eq:I8}
I_8&=\lim_{R\to\infty}[\int_{\partial B_R}\frac{\partial(d-w_0)}{\partial n}\cdot(d-w_0)|d-w_0|^{p-2}d\sigma\\
&-\int_{B_R}|\nabla(d-w_0)|^2|d-w_0|^{p-2}dx\notag\\
&-(p-2)\int_{B_R}|\nabla(d-w_0)\cdot(d-w_0)|^2|d-w_0|^{p-4}dx]\notag\\
&\equiv \lim_{R\to\infty}(K_1-K_2-K_3)\notag.
\end{align}
The boundary term $K_1$ is estimated as
\begin{align}\label{I8boundary}
K_1&=\int_{\partial B_R}\frac{\partial(d-w_0)}{\partial n}\cdot(d-w_0)|d-w_0|^{p-2}d\sigma\\
&\leq\int_{\partial B_R}|\nabla(d-w_0)||d-w_0|^{p-1}d\sigma\notag\\
&\leq\int_{\partial B_R}|\nabla(d-w_0)|d\sigma\notag
\end{align}
where we used  that $|d-w_0|^{p-1}\leq C$ for any $p\geq 1$.\\
In Proposition (\ref{PPGN}), let $k=1$, $m=2$, $r=1$, $p=2$ and $q=1$.  For $a=\frac{2}{7}$  the  inequality (\ref{GN}) yields
\begin{align}\notag
&\int_{\R^3}|\nabla(d(t)-w_0)|dx\\
&\leq C\left(\int_{\R^3}|\Delta(d(t)-w_0)|dx\right)^{1/7}\left(\int_{\R^3}|d(t)-w_0|dx\right)^{5/7}\notag\\
&\leq C_0 (C(t))^{5/7}\notag,
\end{align}
where we used Ladyzhenskaya estimate (\ref{ex:lady}) and the estimate (\ref{ineq:d-w0L1time}), $C_0$ is a constant depending on initial data and $C(t)$ is the time dependent function in (\ref{ineq:d-w0L1time}). Thus, for any fixed $t>0$,
\bg\notag
\int_{\partial B_{R_i}}|\nabla(d-w_0)|d\sigma\to 0
\ed
for a sequence $R_i\to\infty$. Hence, from (\ref{I8boundary}), we have
\bg\label{limit:K1}
\lim_{R\to\infty}K_1=0.
\ed
Since
$$
-|\nabla(d-w_0)|^2|d-w_0|^2\leq -|\nabla(d-w_0)\cdot(d-w_0)|^2,
$$
we have
$$
-|\nabla(d-w_0)|^2|d-w_0|^{p-2}\leq -|\nabla(d-w_0)\cdot(d-w_0)|^2|d-w_0|^{p-4},
$$
which implies from (\ref{eq:I8}) that
\begin{align}\label{est:K2K3}
-K_2-K_3&\leq (p-1)\int_{B_R}-|\nabla(d-w_0)\cdot(d-w_0)|^2|d-w_0|^{p-4}dx\\
&=-\frac{4(p-1)}{p^2}\int_{B_R}|\nabla|d-w_0|^{\frac{p}{2}}|^2dx\notag.
\end{align}
Combining (\ref{eq:I8}), (\ref{limit:K1}) and (\ref{est:K2K3}) yields
\bg\label{limit:I8}
I_8\leq-\frac{4(p-1)}{p^2}\int_{\R^3}|\nabla|d-w_0|^{\frac{p}{2}}|^2dx
\ed
for any $p\geq 1$. Slightly modifying  the process to  estimate $I_6$ gives
\bg\label{eq:I9}
I_9\leq \frac{1}{8\eta^2}\int_{\R^3}|d-w_0|^{p+2}dx.
\ed
Combining  (\ref{eq:Lpd-w0}) with  inequalities (\ref{limit:I7}), (\ref{limit:I8}), and (\ref{eq:I9}) yields
\begin{align}\label{lp:d-w0}
&\frac{1}{p}\frac{d}{dt}\int_{\R^3}|d-w_0|^pdx+(p-1)\int_{\R^3}|\nabla(d-w_0)|^2|d-w_0|^{p-2}dx\\
&\leq C\int_{\R^3}|d-w_0|^{p+2}dx\notag.
\end{align}

Denote $v=|d-w_0|^{p/2}$. Then (\ref{lp:d-w0}) can be rewritten as
\bg\label{lp:v}
\frac{1}{p}\frac{d}{dt}\int_{\R^3}|v|^2dx+\frac{4(p-1)}{p^2}\int_{\R^3}|\nabla v|^2dx\leq C\int_{\R^3}|v|^2|d-w_0|^2dx.
\ed
Applying H$\ddot{o}$lder inequality and Sobolev inequality to the right hand side yields
\begin{align}\notag
\int_{\R^3}|v|^2|d-w_0|^2dx&\leq \left(\int_{\R^3}|v|^6dx\right)^{1/3}\left(\int_{\R^3}|d-w_0|^3dx\right)^{2/3}\\
&\leq C\int_{\R^3}|\nabla v|^2dx\left(\int_{\R^3}|d-w_0|^2dx\right)^{2/3}\|d-w_0\|_{L^\infty(\R^3)}^{2/3}\notag
\end{align}
From Lemma \ref{Le:dl2} and Lemma \ref{Le:linfty}, when $p>1$, we can choose  initial data small enough so that
\bg\label{bad}
C\left(\int_{\R^3}|d-w_0|^2dx\right)^{2/3}\|d-w_0\|_{L^\infty(\R^3)}^{2/3}\leq \frac{2(p-1)}{p^2}.
\ed
It follows from (\ref{lp:v}) that
\bg\label{lp:v2}
\frac{1}{p}\frac{d}{dt}\int_{\R^3}|v|^2dx+\frac{2(p-1)}{p^2}\int_{\R^3}|\nabla v|^2dx\leq 0.
\ed
The inequality (\ref{d-w0lp}), for $p>1$, can be obtained  integrating (\ref{lp:v2}) over  $[0,T]$ .\\

\noindent When $p=1$, we have from (\ref{lp:d-w0})
\begin{align}\notag
\frac{d}{dt}\int_{\R^3}|d-w_0|dx&\leq C\int_{\R^3}|d-w_0|^3dx\\
&\leq C\int_{\R^3}|d-w_0|^2dx\|d-w_0\|_{L^\infty(\R^3)}\notag\\
&\leq C\int_{\R^3}|\nabla d|^2+|\Delta d|^2dx\notag
\end{align}
where we used Lemma \ref{Le:dl2} and Lemma \ref{Le:linfty}. Integrating  over time $[0,T]$ yields
\begin{align}\notag
&\int_{\R^3}|d-w_0|dx\\
&\leq \int_{\R^3}|d_0-w_0|dx+C\int_0^T\int_{\R^3}|\nabla d|^2+|\Delta d|^2dxdt \leq C\notag
\end{align}
 where we used the energy estimate (\ref{basic}). This  completes the proof of the lemma.
\cbdu

\bigskip

\subsection{Decay of $d-w_0$ and $\nabla(d-w_0)$}

We first  establish that $d-w_0$ decays in $L^p(\R^3)$, for $p>1$ at the rate $(1+t)^{-\frac{3}{2}(1-\frac{1}{p})}$.

\begin{Theorem}\label{Thm:d-w0decay}
Let $d$ be the solution obtained in Theorem \ref{exist}. Assume $d-w_0\in L^p(\R^3)$, $p\geq 1$. Assume the initial data satisfies the conditions of  Lemma \ref{Le:dlp}. Then for any $1<p<\infty$, $t>0$
\bg\label{ineq:d-w0decay}
\|d(\cdot,t)-w_0\|_{L^p(\R^3)}\leq C(1+t)^{-\frac{3}{2}(1-\frac{1}{p})},
\ed
where the constant $C$ depends on $\lambda_p$ as defined in Lemma \ref{Le:dlp}.
\end{Theorem}
\pf
Note that as $p \to \infty$ the constants   $\lambda_p$ in Lemma \ref{Le:dlp} will  tend to zero, hence we cannot pass to the limit  as $p\to \infty$. Therefore this result does not give the  decay for the $L^{\infty}$ norm.\\
We proceed by induction for $k$  with $p =2^k$. The other powers $p$ follow by interpolation.
When $k=0$ the theorem follows by Lemma \ref{Le:dlp}. Suppose it holds for $s=k$, then we have
\bg\label{ineq:d-w0decay}
\|d(\cdot,t)-w_0\|_{L^{2^k}(\R^3)}\leq C(1+t)^{-\frac{3}{2}(1-\frac{1}{2^k})},
\ed
Let $v =|d(\cdot,t)-w_0|^{2^k}$.

Recall the inequality (\ref{lp:v2}) ( which holds provided the data satisfies (\ref{bad}))
 \bg\label{cosa}
\frac{d}{dt}\int_{\R^3}|v|^2dx+C_p\frac{2(p-1)}{p}\int_{\R^3}|\nabla v|^2dx\leq 0,
\ed
By Gagliardo-Nirenberg we have

\[\int_{\R^3}|v|^2dx \leq C \left(\int |\nabla v|^2 dx \right)^{3/5}\left( \int|v| dx  \right)^{4/5}\]
Hence using the inductive hypothesis on the last  integral on the right hand side  we have
\[ -\int |\nabla v|^2 dx  \leq - C\left( \int_{\R^3}|v|^2dx\right)^{5/3} ( 1+t)^{2(2^k-1)}\]
Combining the last inequality with  (\ref{cosa}) yields
\[\frac{d}{dt} \frac{\int_{\R^3}|v|^2dx}{\left( \int_{\R^3}|v|^2dx\right)^{5/3}} \leq -( 1+t)^{2(2^k-1)}\]
Integrating  and reordering terms yields

\[ \int_{\R^3}|v|^2dx \leq \left[ \frac{v_0}{ (1+(1+t) )^{2(2^k-1)+1}}\right]^{3/2}\]
Since $\frac{3}{2}(2(2^k-1)+1) =\frac{3}{2}(2^{k+1} -1)$ the induction step is obtained, establishing the conclusion of the theorem.

\cbdu

As a consequence of the last theorem, we derive the decay of $\nabla(d-w_0)$.
\begin{Corollary}\label{dh1decay}Let $d$ be the solution to system (\ref{LCD}) obtained in Theorem \ref{exist}. Then
\bg\notag
\|\nabla(d-w_0)\|_{L^2}^2\leq C(1+t)^{-\frac{3}{4}},
\ed
where $C$ depends on initial data.
\end{Corollary}
\pf
Take $p=2$ in Theorem \ref{Thm:d-w0decay},
\bg\label{ineq:d-w0decayl2}
\|d(\cdot,t)-w_0\|_{L^2(\R^3)}\leq C\|d_0-w_0\|_{L^2(\R^3)}(1+t)^{-\frac{3}{4}}.
\ed
Gagliardo-Nirenberg inequality (\ref{GN}) yields
\begin{align}\notag
\int_{\R^3}|\nabla(d-w_0)|^2dx&\leq C\left (\int_{\R^3}|d-w_0|^2dx\right )^{\frac{1}{2}}\left (\int_{\R^3}|\Delta(d-w_0)|^2dx\right )^{\frac{1}{2}}\\
&\leq C\left (\int_{\R^3}|d-w_0|^2dx\right )^{\frac{1}{2}}\notag\\
&\leq C(1+t)^{-\frac{3}{4}}\notag,
\end{align}
where in the last two steps we used Ladyzhenskaya energy estimate (\ref{lady}) and (\ref{ineq:d-w0decayl2}), respectively. The constant $C$ depends on initial data. It completes the proof.
\cbdu

\bigskip
\section{Decay of Velocity}

 An application of the Fourier Splitting method  \cite{Sch} is used to  establish $L^2$ decay of velocity $u$.

\begin{Theorem}\label{ul2decay}
Let $u$ be the solution obtained in Theorem \ref{exist}. If additionally $u_0\in L^1(\R^3)$, then
\bg\notag
\|u(\cdot,t)\|_{L^2}^2\leq C(1+t)^{-\frac{1}{2}},
\ed
where $C$ depends on initial data, the $L^1$ and $L^2$ norm of $u_0$.
\end{Theorem}
\pf
Multiplying the Navier-Stokes equation in system (\ref{LCD}) by $u$ and integrating by parts yields
\bg\label{NSEen}
\frac{d}{dt}\int_{\R^3}|u|^2+2\int_{\R^3}|\nabla u|^2dx=2\int_{\R^3}\nabla u(\nabla d\otimes\nabla d)dx.
\ed
H\"older and Cauchy Schwartz inequalities  yield
\bg\notag
2\int_{\R^3}\nabla u(\nabla d\otimes\nabla d)dx\leq \int_{\R^3}|\nabla u|^2dx+C\int_{\R^3}|\nabla d\otimes\nabla d|^2dx.
\ed
Thus, we derive from (\ref{NSEen})
\bg\notag
\frac{d}{dt}\int_{\R^3}|u|^2+\int_{\R^3}|\nabla u|^2dx\leq C\int_{\R^3}|\nabla d\otimes\nabla d|^2dx.
\ed
The right hand side of above inequality can be estimated as
\begin{align}\notag
\int_{\R^3}|\nabla d\otimes\nabla d|^2dx&=\int_{\R^3}(\nabla d\otimes\nabla d)(\nabla d\otimes\nabla d)dx\\
&=-3\int_{\R^3}(d-w_0)\;\Delta d\;\nabla d\otimes\nabla d\notag\\
&\leq \frac{1}{2}\int_{\R^3}|\nabla d\otimes\nabla d|^2dx+C\int_{\R^3}|d-w_0|^2|\Delta d|^2dx\notag,
\end{align}
from which it follows
\begin{align}\notag
\int_{\R^3}|\nabla d\otimes\nabla d|^2dx&\leq C\int_{\R^3}|d-w_0|^2|\Delta d|^2dx\\
&\leq C\left(\int_{\R^3}|d-w_0|^pdx\right)^{\frac{2}{p}}\left(\int_{\R^3}|\Delta d|^{\frac{2p}{p-2}}dx\right)^{\frac{p-2}{p}}\notag\\
&=C\|d-w_0\|_{L^p(\R^3)}^2\left(\int_{\R^3}|\Delta d|^{2+\frac{4}{p-2}}dx\right)^{\frac{p-2}{p}}\notag\\
&\leq C\|d-w_0\|_{L^p(\R^3)}^2\notag,
\end{align}
for $p\geq 2$, where the last step followed from Ladyzhenskaya estimate (\ref{lady}) and the fact that $\|\Delta d\|_{L^\infty(\R^3\times[0,T])}$ is bounded since $d$ is regular in the sense stated in Theorem \ref{exist}. Thus, it follows from Theorem \ref{Thm:d-w0decay} that
\bg\notag
\int_{\R^3}|\nabla d\otimes\nabla d|^2dx\leq C(1+t)^{-3(1-\frac{1}{p})},
\ed
for any $p\geq 2$. Therefore,
\bg\label{NSEen1}
\frac{d}{dt}\int_{\R^3}|u|^2dx+\int_{\R^3}|\nabla u|^2dx\leq C(1+t)^{-3(1-\frac{1}{p})}.
\ed
Applying Plancherel's theorem to (\ref{NSEen1}) gives
\bg\label{NSEF}
\frac{d}{dt}\int_{\R^3}|\hat{u}|^2d\xi+\int_{\R^3}|\xi|^2|\hat{u}|^2d\xi\leq C(1+t)^{-3(1-\frac{1}{p})}.
\ed
The idea is to decompose the frequency domain $\R^3$ in integral $\int_{\R^3}|\xi|^2|\hat{u}|^2d\xi$ into two time-dependent subdomains. The time dependent subdomains are a $3$-dimensional sphere, $S(t)$, centered at the origin with an appropriate time dependent radius and its complement. For this we rewrite (\ref{NSEF})  as
\bg\notag
\frac{d}{dt}\int_{\R^3}|\hat{u}|^2d\xi\leq-\int_{S(t)^c}|\xi|^2|\hat{u}|^2d\xi-\int_{S(t)}|\xi|^2|\hat{u}|^2d\xi
+C(1+t)^{-3(1-\frac{1}{p})},
\ed
where $S(t)$ is the ball
\[S(t) = \{\xi \in \R^3: |\xi|\leq r(t)=(\frac{k}{1+t})^{1/2}\}\]
for a certain $k$, which will be determined below. Hence
\begin{align}\notag
\frac{d}{dt}\int_{\R^3}|\hat{u}|^2d\xi&\leq-\frac{k}{1+t}\int_{S(t)^c}|\hat{u}|^2d\xi-\int_{S(t)}|\xi|^2|\hat{u}|^2d\xi
+C(1+t)^{-3(1-\frac{1}{p})}\\
&=-\frac{k}{1+t}\int_{\R^3}|\hat{u}|^2d\xi+\int_{S(t)}(\frac{k}{1+t}-|\xi|^2)|\hat{u}|^2d\xi
+C(1+t)^{-3(1-\frac{1}{p})}\notag
\end{align}
and
\bg\label{NSEF1}
\frac{d}{dt}\int_{\R^3}|\hat{u}|^2d\xi+\frac{k}{1+t}\int_{\R^3}|\hat{u}|^2d\xi\leq
\frac{k}{1+t}\int_{S(t)}|\hat{u}|^2d\xi+C(1+t)^{-3(1-\frac{1}{p})}.
\ed
The following estimate, which will be established later, is needed
\bg\label{Fu}
|\hat{u}(\xi,t)|\leq C|\xi|^{-1}
\ed
for $\xi\in S(t)$, where $C$ is a constant only depending on the initial data. Combining the inequalities (\ref{NSEF1}) and (\ref{Fu}) yields
\bg\notag
\frac{d}{dt}\int_{\R^3}|\hat{u}|^2d\xi+\frac{k}{1+t}\int_{\R^3}|\hat{u}|^2d\xi\leq
\frac{C}{1+t}\int_{S(t)}|\xi|^{-2}d\xi+C(1+t)^{-3(1-\frac{1}{p})}.
\ed
Multiplying by the integrating factor $(1+t)^k$ yields
\bg\notag
\frac{d}{dt}\left[(1+t)^k\int_{\R^3}|\hat{u}|^2d\xi\right]\leq C(1+t)^{k-1}\int_{S(t)}|\xi|^{-2}d\xi+C(1+t)^{k-3(1-\frac{1}{p})}.
\ed
Since $p \geq 2$ and
\bg\notag
\int_{S(t)}|\xi|^{-2}d\xi\leq C\int_0^{r(t)}r^2r^{-2}dr\leq C(1+t)^{-1/2},
\ed
it follows that
\bg\notag
\frac{d}{dt}\left[(1+t)^k\int_{\R^3}|\hat{u}|^2d\xi\right]\leq C(1+t)^{k-\frac{3}{2}}.
\ed
Integrating in time yields
\bg\notag
(1+t)^{k}\int_{\R^3}|\hat{u}|^2d\xi\leq\int_{\R^3}|\hat{u}(\xi,0)|^2d\xi+C[(1+t)^{k-\frac{1}{2}}-1].
\ed
Thus,
\bg\notag
\int_{\R^3}|\hat{u}|^2d\xi\leq(1+t)^{-k}\int_{\R^3}|\hat{u}(\xi,0)|^2d\xi+C[(1+t)^{-\frac{1}{2}}-(1+t)^{-k}].
\ed
Since $u_0\in L^2$, it follows that $\hat u(0)\in L^2$ from Plancherel's theorem. Hence
\bg\notag
\int_{\R^3}|\hat{u}|^2d\xi\leq C(1+t)^{-\frac{1}{2}},
\ed
Hence
\bg\notag
\int_{\R^3}|u|^2dx\leq C(1+t)^{-\frac{1}{2}}.
\ed

To complete the proof we need to establish the inequality (\ref{Fu}). Taking the Fourier transform of Navier-Stokes equation in system (\ref{LCD}) yields
\bg\label{FNSE}
\hat u_t+|\xi|^2\hat u=G(\xi,t)
\ed
where
\bg\notag
G(\xi,t)=-\mathcal F(u\cdot\nabla u)-\mathcal F(\nabla p)-\mathcal F(\nabla\cdot(\nabla d\otimes\nabla d)),
\ed
and $\mathcal F$ indicates the Fourier transform. Multiplying (\ref{FNSE}) by the integrating factor $e^{|\xi|^2t}$ yields
\bg\notag
\frac{d}{dt}[e^{|\xi|^2t}\hat u]=e^{|\xi|^2t}G(\xi,t).
\ed
Integrating in time gives
\bg\label{Fueq}
\hat u(\xi,t)=e^{-|\xi|^2t}\hat u_0+\int_0^te^{-|\xi|^2(t-s)}G(\xi,s)ds.
\ed
We assume for the moment the following auxiliary estimate, which  we will prove below,
\bg\label{G}
|G(\xi,t)|\leq C|\xi|.
\ed
Combining (\ref{Fueq}) and (\ref{G}) yields
\bg\label{Fuineq}
|\hat u(\xi,t)|\leq e^{-|\xi|^2t}|\hat u_0|+\int_0^te^{-|\xi|^2(t-s)}|\xi|ds.
\ed
Since $u_0\in L^1$, we have $|\hat u_0|\leq C$ for all $\xi$ and some constant $C$.Performing integration in (\ref{Fuineq}) gives
\bg\notag
|\hat u(\xi,t)|\leq Ce^{-|\xi|^2t}+\frac{C}{|\xi|}(1-e^{-|\xi|^2t})\leq C|\xi|^{-1}
\ed
for $\xi\in S(t)$. To finish the proof we need to establish (\ref{G}). For this purpose we analyze each term in $G(\xi,t)$ separately. We have
\bg\notag
|\mathcal F(u\cdot\nabla u)|=|\mathcal F(\nabla\cdot(u\otimes u))|\leq \sum_{i,j}\int_{\R^3}|u^iu^j||\xi_j|dx.
\ed
Since $u\in L^\infty(L^2)$ by the basic energy estimate, we have
\bg\notag
|\mathcal F(u\cdot\nabla u)|\leq C|\xi|.
\ed
 By the basic energy inequalities (\ref{basic}, \ref{lady})  we have  $\nabla d\in L^\infty(L^2)$ proceeding  similarly as for  the last inequality we have
\bg\notag
|\mathcal F(\nabla\cdot(\nabla d\otimes\nabla d))|\leq C|\xi|.
\ed
Taking divergence of Navier-Stokes equation in system (\ref{LCD}) gives that
\bg\notag
\Delta p=-\sum_{i,j}\frac{\partial^2}{\partial x_i\partial x_j}(u^iu^j)-\sum_{i,j}\frac{\partial^2}{\partial x_i\partial x_j}(\nabla d^i\nabla d^j).
\ed
Taking the Fourier transform then yields
\bg\notag
|\xi|^2\mathcal F(p)=-\sum_{i,j}\xi_i\xi_j\mathcal F(u^iu^j)-\sum_{i,j}\xi_i\xi_j\mathcal F(\nabla d^i\nabla d^j).
\ed
Since $\mathcal F(u^iu^j)\in L^\infty$ and $\mathcal F(\nabla d^i\nabla d^j)\in L^\infty$, it follows that
\bg\notag
\mathcal F(p)\leq C,
\ed
and thus $\mathcal F(\nabla p)\leq C|\xi|$. It completes the proof of (\ref{G}) and hence completes the proof of theorem.
\cbdu

 Combining Theorem \ref{Thm:d-w0decay}, Corollary \ref{dh1decay} and Theorem \ref{ul2decay} yields the proof of the main Theorem \ref{Mthm}.

 \cbdu

\bigskip

\begin{Remark}
The decay rate for the velocity $u$ in $L^2$ obtained in \cite{Wu}, for the bounded domain case, is $(1+t)^{-\frac{\theta}{1-2\theta}}$  where $\theta\in(0,\frac{1}{2})$. When $\theta$ is close to $0$,  then $-\frac{\theta}{1-2\theta}$ would be very small, meaning the decay is very slow. In this paper, we obtained the decay rate for velocity $u$ in $L^2$ with $(1+t)^{-\frac{1}{4}}$, a fixed constant algebraic rate. The advantage comes from the fact that we work on the whole space $\R^3$ where we can apply the Fourier splitting method.
\end{Remark}

\begin{Remark}
 It was pointed out in the first section that there is an essential difficulty to apply \L ojasiewicz-Simon approach in whole space $\R^3$. However, in weighted Sobolev spaces of $\R^3$, the compactness is recovered. Thus, we expect  there is  hope to construct certain \L ojasiewicz-Simon type inequality in weighted Sobolev spaces and proceed  with the method in \cite{Wu} to derive the decay of solutions to the LCD system in weighted Sobolev spaces.
\end{Remark}

\bigskip

\appendix

\section {Existence of Classical Solutions in $\R^3$} \label{Ap}

In this section we sketch a brief proof of the existence Theorem \ref{exist},  Section $1$.
As mentioned in the introduction, for  bounded domains in $\R^3$, the existence of global regular solutions to the flow of nematic liquid crystals with constant density has been established in \cite{LL} provided the viscosity is large enough. The existence of global regular solutions to the flow of nematic liquid crystals with non-constant density has been established in \cite{DQS} provided the initial data is small enough. In both of the above papers, a Ladyzhenskaya energy estimate (higher order derivative estimate) was derived  and hence a relatively  standard bootstrapping argument  yielded a regular solution.\\

The proof of Theorem \ref{exist} will be given through four steps. In the first step, on a sequence of balls $B_{R_n}$ with radius $R_n$, centered at the origin, we obtain the existence of a Galerkin approximated solution $(u^{n,m},d^{n,m})$ for the system (\ref{LCD}) with modified initial data, for each $m=1,2,3,...$. In the second step, we establish an estimate of $d^{n,m}-w_0$ in $L^1(B_{R_n})$ for any fixed time $t>0$. In the third step, we take the limit $m\to\infty$. In the forth step, we take the limit $R_n\to\infty$. In fact, we are able to show  that all the estimates in Ladyzhenskaya energy method in step one are independent of  the domain size. Thus we can take a subsequence of solutions on balls $B_{R_n}$ which converge to a limit in $\R^3$ when $R_n$ goes to infinity.\\

\begin{Lemma}\label{Le:initial}
Assume $u_0\in H^1(\R^3)$ and $\tilde d_0\equiv d_0(x)-w_0\in H^2(\R^3)\cap L^1(\R^3)$ with $|w_0|=1$. There exists a sequence of functions $\left\{(u_0^n, d_0^n)\right\}_{n=1}^\infty$ and a sequence of real numbers $\left\{R_n\right\}_{n=1}^\infty$ with $R_n\to\infty$ as $n\to\infty$ such that,
\bg\label{limit:uinitial}
u_0^n\in H_0^1(B_{R_n}),  \text{ with } u_0^n\to u_0 \text{ in } H^1(\R^3) \text{ as } n\to\infty
\ed
\bg\label{limit:dinitial}
\tilde d_0^n\in H_0^2(B_{R_n}),  \text{ with } \tilde d_0^n\to \tilde d_0 \text{ in } H^2(\R^3) \text{ as } n\to\infty
\ed
where $\tilde d_0^n=d_0^n-w_0$. Moreover,
\bg\label{app:initialL1}
\|\tilde d_0^n\|_{L^1(B_{R_n})}\leq C\|\tilde d_0\|_{L^1(\R^3)},
\ed
and $|d_0^n|\leq 1$.
\end{Lemma}
\pf
Such a sequence of functions can be constructed easily as follows. Let $\zeta_n$ be a sequence of  smooth functions such that
\begin{equation}\notag
\zeta_n(x) =
\begin{cases}
  1,  & \mbox{ in } B_{R_n} \\
  0, & \mbox{ in } B_{2R_n}
\end{cases}
\end{equation}
and $|\zeta_n(x)|\leq 1$ for all $x\in\R^3$. \\
Define $u_0^n=\zeta_nu_0$ and $\tilde d_0^n=\zeta_n\tilde d_0$. Let $d_0^n=\tilde d_0^n+w_0$. These $\zeta_n$ can be chosen such that
(\ref{limit:uinitial}), (\ref{limit:dinitial}) and (\ref{app:initialL1}) are satisfied. In addition, we have
\bg\notag
|d_0^n|=|\zeta_nd_0+(1-\zeta_n)w_0|\leq \zeta_n+(1-\zeta_n)=1.
\ed
\cbdu

 For the sequel we assume as   initial conditions
\bg\label{initial-A}
\left\{\aligned u(x,0) & = u_0^n(x), \\  d(x, 0) & =d_0^n(x), \endaligned\right., \text{ in $B_{R_n}\times [0, T)$},
\ed
where $u_0^n$ and $d_0^n$ are as obtained in  Lemma \ref{Le:initial}, and the boundary conditions
\bg\label{bd-A}
\left\{\aligned u(x,t) & = 0, \\  d(x, t)-w_0 & =0, \endaligned\right., \text{ on $\partial B_{R_n}\times [0, T)$}.
\ed

We have the following existence of solutions to system (\ref{LCD}) with (\ref{initial-A}) and (\ref{bd-A}) on each ball $B_{R_n}$.
\begin{Theorem}\label{existbdGK} Let $B_{R_n}$ be the ball centered at origin with radius $R_n$ in $\R^3$. Assume that $u_0\in H^1(\R^3)$ and $d_0-w_0\in H^{2}(\R^3)$.
The system (\ref{LCD}) with initial and boundary conditions (\ref{initial-A}) and (\ref{bd-A}) has a smooth solution $(u^{n,m},p^{n,m},d^{n,m})$ for each $m=1,2,3, ...$ satisfying, for any $T>0$
\bg\notag
u^{n,m}\in L^2(0,T;H_0^1(B_{R_n}))\cap L^\infty(0,T;L^2(B_{R_n}))
\ed
\bg\notag
\tilde d^{n,m}\in L^2(0,T;H_0^2(B_{R_n}))\cap L^\infty(0,T;H_0^1(B_{R_n}))
\ed
with $\tilde d^{n,m}=d^{n,m}-w_0$ and $|d^{n,m}|\leq 1$. Moreover, it satisfies the energy inequality
\begin{align}\label{energy-nm}
&\frac{d}{dt}\int_{B_{R_n}}\frac{1}{2}|u^{n,m}|^2+\frac{1}{2}|\nabla d^{n,m}|^2+F(d^{n,m})dx\\
&+\int_{B_{R_n}}|\nabla u^{n,m}|^2+|\Delta d^{n,m}-f(d^{n,m})|^2dx\leq 0.\notag
\end{align}
\end{Theorem}
\pf
The existence proof is obtained through the standard Galerkin approximation method, see \cite{LL} and \cite{DQS}. We only need to give a brief explanation on the claim $|d^{n,m}|\leq 1$ by applying a maximum principle argument. Notice that the approximated initial data $d_0^n$ satisfies $|d_0^n(x)|\leq 1$ for all $x\in\R^3$. Suppose there exists a point $(x_0,t_0)$ in the interior of the domain $B_{R_n}\times[0,T)$, such that $|d^{n,m}|^2$ attains a maximum value at this point. Multiplying the equation
$$
d_t^{n,m}+u^{n,m}\cdot\nabla d^{n,m}=\Delta d^{n,m}-\frac{1}{\eta^2}(|d^{n,m}|^2-1)d^{n,m}
$$
by $d^{n,m}$ yields
\begin{align}\label{normd}
&\frac{d}{dt}|d^{n,m}|^2+u^{n,m}\cdot\nabla|d^{n,m}|^2\cdot d^{n,m}\\
&=\triangle |d^{n,m}|^2
-2|\nabla d^{n,m}|^2-\frac{2}{\eta^2}(|d^{n,m}|^2-1)|d^{n,m}|^2\notag.
\end{align}
At the maximum point $(x_0,t_0)$, we have $\frac{d}{dt}|d^{n,m}|^2=\nabla|d^{n,m}|^2=0$ and $\triangle |d^{n,m}|^2\leq 0$. Thus, it follows from the equation (\ref{normd}) that, at the point $(x_0,t_0)$
$$
(|d^{n,m}|^2-1)|d^{n,m}|^2\leq 0.
$$
This insures that $|d^{n,m}|\leq 1$ at any  interior maximum point $(x_0,t_0)$. Therefore, $|d^{n,m}|\leq 1$ for all points in $B_{R_n}\times[0,T)$.
\cbdu

For the solution $(u^{n,m}, d^{n,m})$ obtained in the above theorem on $B_{R_n}\times[0,T)$, we define energy quantity
\bg\notag
\Phi_{n,m}^2(t)=\|\nabla u^{n,m}\|_{L^2(B_{R_n})}^2+\|\triangle (d^{n,m}-w_0)\|_{L^2(B_{R_n})}^2.
\ed
With a slight  modification of the proof of Theorem 3.1 in \cite{DQS}, we are able to show  that

\begin{Theorem}\label{Thm:lady-nm}Assume that $u_0\in H^1(B_{R_n})$ and $d_0\in H^{2}(B_{R_n})$, and $\|u_0\|_{ H^1(\R^3)}^2+\|\tilde d_0\|_{H^{2}(\R^3)}^2<\infty$. Let $(u^{n,m}, d^{n,m})$ be solutions obtained in Theorem \ref{existbdGK}.
There is a positive small number $\e_0$ such that if
\begin{equation}\label{smalli}
\|u_0\|_{H^1(\R^3)}^2+\|d_0- w_0\|_{H^2(\R^3)}^2\leq\epsilon_0,
\end{equation}
then
\begin{align}\label{energy:Lady-nm}
&\int_{B_{R_n}}|\nabla u^{n,m}|^2+|\Delta (d^{n,m}-w_0)|^2dx\\
&+\int_0^T\int_{B_{R_n}}|\Delta u^{n,m}|^2+|\nabla\Delta (d^{n,m}-w_0)|^2dx\notag\\
&\leq C(\|u_0\|_{H^1(B_{R_n})}^2+\|d_0- w_0\|_{H^2(B_{R_n})}^2)\notag,
\end{align}
for any $T>0$, where the constant $C$ is independent of domain size $R_n$ and $m$.
\end{Theorem}

There is no need to prove the theorem except that we need a brief explanation on the last claim that constant $C$ is independent of $R_n$. In the proof of Ladyzhenskaya energy estimate  in \cite{DQS}, we only  use the Gagliardo-Nirenberg interpolation inequalities and standard elliptic inequalities. That is we use
\bg\notag
\|u^{n,m}\|_{L^4}^4\leq C\|u^{n,m}\|_{L^2}\|\nabla u^{n,m}\|_{L^2}^3
\ed
\bg\notag
\|\nabla (d^{n,m}-w_0)\|_{L^4}^4\leq C\|\nabla (d^{n,m}-w_0)\|_{L^2}\|\Delta (d^{n,m}-w_0)\|_{L^2}^3
\ed
and the  elliptic estimate
\bg\notag
\|D^2u^{n,m}\|_{L^2}\leq C\|\Delta u^{n,m}\|_{L^2}
\ed
for $u^{n,m}$ and $d^{n,m}-w_0$ vanishing on the boundary. In the above inequalities, the various constants $C$ are independent of  the size of the domain. Thus the constant $C$ in (\ref{energy:Lady-nm}) is independent of $R_n$.\\

For  the second step, we derive a time dependent estimate of $d^{n,m}(t)-w_0$ in $L^1(B_{R_n})$.

\begin{Lemma}\label{Le:dL1-nm}
Let $d^{n,m}$ be the solution obtained in Theorem \ref{existbdGK}. In addition, assume $d_0-w_0\in L^1(\R^3)$. Then
\bg\label{ineq:d-w0L1}
\int_{B_{R_n}}|d^{n,m}(t)-w_0|dx\leq (C_0t+\int_{\R^3}|d_0-w_0|dx)e^{Ct}
\ed
where the constant $C_0$ only depends on initial data and constant $C$ only depends on $\eta$.
\end{Lemma}
\pf
By the second equation in (\ref{LCD}) we have,
\begin{align}\label{d-w0:L1}
&\frac{d}{dt}\int_{B_{R_n}}|d^{n,m}(t)-w_0|dx\\
&=\int_{B_{R_n}}\frac{(d^{n,m}(t)-w_0)\cdot d^{n,m}_t}{|d^{n,m}(t)-w_0|}dx\notag\\
&=\int_{B_{R_n}}\frac{(d^{n,m}(t)-w_0)\cdot \Delta d^{n,m}(t)}{|d^{n,m}(t)-w_0|}dx-\int_{B_{R_n}}\frac{(d^{n,m}(t)-w_0)\cdot[u^{n,m}\cdot\nabla d^{n,m}(t)]}{|d^{n,m}(t)-w_0|}dx\notag\\
&-\int_{B_{R_n}}\frac{(d^{n,m}(t)-w_0)\cdot f(d^{n,m}(t))}{|d^{n,m}(t)-w_0|}dx\notag\\
&\equiv I_1+I_2+I_3.\notag
\end{align}
There is no need to worry about the singular points of $(d^{n,m}-w_0)^{-1}$ in the above equation, since each term on the right hand side contains $\frac{d^{n,m}-w_0}{|d^{n,m}-w_0|}$. We deal with the three terms $I_1$, $I_2$ and $I_3$ in the following way.\\
Since $d^{n,m}(t)-w_0=0$ on the boundary $\partial B_{R_n}$. Replacing $d_t$ by its value in the direction equation and  integration by parts yields
\begin{align}\label{estimate:I}
I_1&=-\int_{B_{R_n}}\frac{|\nabla(d^{n,m}(t)-w_0)|^2}{|d^{n,m}(t)-w_0|}dx \\
&+\int_{B_{R_n}}|(d^{n,m}(t)-w_0)\cdot \nabla (d^{n,m}(t)-w_0)|^2|d^{n,m}(t)-w_0|^{-3}dx\notag\\
&\leq 0,\notag
\end{align}
since $|(d^{n,m}(t)-w_0)\cdot \nabla (d^{n,m}(t)-w_0)|^2\leq |(d^{n,m}(t)-w_0)|^2|\nabla (d^{n,m}(t)-w_0)|^2$.\\
By  H\"o lder inequality we have
\begin{align}\label{estimate:II}
|I_2|&\leq \int_{B_{R_n}}|u^{n,m}||\nabla(d^{n,m}-w_0)|dx\\
&\leq C\left(\int_{B_{R_n}}|u^{n,m}|^2dx\right)^{1/2}\left(\int_{B_{R_n}}|\nabla d^{n,m}|^2dx\right)^{1/2}\leq C_0\notag,
\end{align}
where we used the energy estimate (\ref{energy-nm}), and the constant $C_0$ only depends on the initial data.\\
Recall that by definition $f(d^{n,m})=\frac{1}{\eta^2}(|d^{n,m}|^2-1)d^{n,m}$, and  $|d^{n,m}|\leq 1$ from Theorem \ref{existbdGK} and $|w_0|=1$, hence
\begin{align}\label{estimate:III}
|I_3|&\leq \frac{1}{\eta^2}\int_{B_{R_n}}|f(d^{n,m})|dx\\
&\leq C\int_{B_{R_n}}|d^{n,m}-w_0||d^{n,m}+w_0||d^{n,m}|dx\leq C\int_{B_{R_n}}|d^{n,m}-w_0|dx\notag,
\end{align}
where the constant $C$ depends on $\eta$.
Combining the inequalities (\ref{d-w0:L1}), (\ref{estimate:I}), (\ref{estimate:II}) and (\ref{estimate:III}) yields
\begin{equation}\notag
\frac{d}{dt}\int_{B_{R_n}}|d^{n,m}(t)-w_0|dx\leq C\int_{B_{R_n}}|d^{n,m}(t)-w_0|dx+C_0.
\end{equation}
Integrating  over $[0,t]$,  and  Gronwall's inequality (see \cite{Ev}) gives,
\bg\notag
\int_{B_{R_n}}|d^{n,m}(t)-w_0|dx\leq (C_0t+\int_{B_{R_n}}|d_0^n-w_0|dx)e^{Ct},
\ed
for any $t>0$.
Inequality (\ref{ineq:d-w0L1}) now follows from the last estimate in Lemma \ref{Le:initial}. This concludes the second step.\\

In the third step, we take the limit of the Galerkin approximating solutions $(u^{n,m},d^{n,m})$ as $m\to\infty$.
By the estimates (\ref{energy-nm}) and (\ref{energy:Lady-nm}), there exists $(u^n,d^n)$ for each $n=1,2,3, ...$ such that, taking subsequence if necessary,
$$
u^{n,m}\rightharpoonup u^n \ \mbox { weakly in } L^2(0,T;H_0^1(B_{R_n})),
$$
$$
u^{n,m}\to u^n \ \mbox { strongly in } L^\infty(0,T;L^2(B_{R_n})),
$$
$$
d^{n,m}\rightharpoonup d^n \ \mbox { weakly in } L^2(0,T;H^2(B_{R_n})),
$$
$$
d^{n,m}\to d^n  \ \mbox { strongly in } L^2(0,T;H^1(B_{R_n})) \text{ with } d^n=w_0 \text{ on } \partial{B_{R_n}}.
$$
It follows easily from the above convergence that $(u^n, d^n)$ is a weak solution to the system \eqref{LCD} with initial condition (\ref{initial-A}) and boundary condition (\ref{bd-A}) on $B_{R_n}\times[0,T)$. Moreover,  the solutions $(u^n,d^n)$ satisfy the basic energy inequality
\begin{align}\label{energy-n}
&\frac{d}{dt}\int_{B_{R_n}}\frac{1}{2}|u^{n}|^2+\frac{1}{2}|\nabla d^{n}|^2+F(d^{n})dx\\
&+\int_{B_{R_n}}|\nabla u^{n}|^2+|\Delta d^{n}-f(d^{n})|^2dx\leq 0,\notag
\end{align}
and the higher order energy inequality
\begin{align}\label{energy:Lady-n}
&\int_{B_{R_n}}|\nabla u^{n}|^2+|\Delta (d^{n}-w_0)|^2dx\\
&+\int_0^t\int_{B_{R_n}}|\Delta u^{n}|^2+|\nabla\Delta (d^{n}-w_0)|^2dx\notag\\
&\leq C(\|u_0\|_{H^1(B_{R_n})}^2+\|d_0- w_0\|_{H^2(B_{R_n})}^2)\notag,
\end{align}
for any $t>0$, where the constant $C$ is independent of domain size $R_n$. In addition, from Lemma \ref{Le:dL1-nm} it follows that $d^n$ satisfies the estimate
\bg\label{ineq:d-w0L1n}
\int_{B_{R_n}}|d^{n}(t)-w_0|dx\leq (C_0t+\int_{\R^3}|d_0-w_0|dx)e^{Ct}
\ed
where the constant $C_0$ only depends on initial data and constant $C$ only depends on $\eta$.\\

In the forth step, we extend the solutions $(u^n,d^n)$ on $B_{R_n}$ to the whole space $\R^3$ by taking limit $R_n\to\infty$. With the estimates (\ref{energy-n}) and (\ref{energy:Lady-n}) we can extract a subsequence $\left\{(u^{1k},d^{1k})\right\}_{k=1}^\infty$ from $(u^n, d^n)$  for $n\geq 1$ such that
\bg\notag
u^{1k}\rightharpoonup u_{(1)} \ \mbox{in } L^2(0,T;H^1(B_{R_1}))
\ed
\bg\notag
u^{1k}\to u_{(1)} \ \mbox{in } L^\infty(0,T;L^2(B_{R_1}))
\ed
\bg\notag
d^{1k}\rightharpoonup d_{(1)} \ \mbox{in } L^2(0,T;H^2(B_{R_1}))
\ed
\bg\notag
d^{1k}\to d_{(1)} \ \mbox{in } L^2(0,T;H^1(B_{R_1}))
\ed
and the limit $(u_{(1)},d_{(1)})$ satisfies system (\ref{LCD}) in distribution sense and the estimates (\ref{energy-n}) (\ref{energy:Lady-n})  on $B_{R_1}\times[0,T)$.

On $B_{R_2}\times[0,T)$, we take subsequence $\left\{(u^{2k},d^{2k})\right\}_{k=1}^\infty$ from $\left\{(u^{1k},d^{1k})\right\}_{k=1}^\infty$ such that $\left\{u^{2k}\right\}_{k=1}^\infty$, $\left\{p^{2k}\right\}_{k=1}^\infty$ and $\left\{d^{2k}\right\}_{k=1}^\infty$ converge to $u_{(2)}$ and $d_{(2)}$ respectively in the same convergence sense as above. And we have that
$$
u_{(2)}|_{B_1}=u_{(1)}, \ \ d_{(2)}|_{B_1}=d_{(1)}.
$$
Repeating  the process on each  $B_{R_n}\times[0,T)$, we can take subsequence $\left\{(u^{nk},d^{nk})\right\}_{k=1}^\infty$ from the sequence $\left\{(u^{(n-1)k},d^{(n-1)k})\right\}_{k=1}^\infty$, such that $\left\{u^{nk}\right\}_{k=1}^\infty$ and $\left\{d^{nk}\right\}_{k=1}^\infty$ converge to $u_{(n)}$ and $d_{(n)}$ respectively. And we have that
$$
u_{(n)}|_{B_{n-1}}=u_{(n-1)}, \ \ d_{(n)}|_{B_{n-1}}=d_{(n-1)}.
$$
Then we take the diagonal sequence $\left\{(u^{kk},d^{kk})\right\}_{k=1}^\infty$ and let $k\to\infty$. This sequence (if necessary, take a subsequence of it) converges to $(u,d)$, in $\mathbb{R}^3\times[0,T)$. The limit $(u,d)$ satisfies the system (\ref{LCD}) in the  sense of distributions  and satisfies the energy estimates
\begin{align}\label{limitbasicenergy}
&\int_{\R^3}|u|^2+|\nabla d|^2+2F(d)dx+2\int_0^T\int_{\R^3}|\nabla u|^2+|\Delta d-f(d)|^2dxdt\\
&\leq \int_{\R^3}|u_0|^2+|\nabla d_0|^2dx\notag
\end{align}
\begin{align}\label{limitLadyenergy}
&\int_{\R^3}|\nabla u|^2+|\Delta d|^2dx+\int_0^T\int_{\R^3}|\Delta u|^2+|\nabla\Delta d|^2dxdt\\
&\leq C(\|u_0\|_{H^1(\R^3)}^2+\|d_0- w_0\|_{H^2(\R^3)}^2).\notag
\end{align}
In addition, the solution $d$ satisfies the estimate
\bg\notag
\int_{\R^3}|d(t)-w_0|dx\leq (C_0t+\int_{\R^3}|d_0-w_0|dx)e^{Ct}
\ed
with constant $C_0$ only depending on initial data and constant $C$ only depending on $\eta$.\\

Estimates (\ref{limitbasicenergy}) and (\ref{limitLadyenergy}) allow us to  apply the ``bootstrapping argument"  as used in \cite{DQS} and \cite{LL},  and  prove that the limit $(u,p,d)$ is a classical solution to system (\ref{LCD}) satisfying the desired estimates in Theorem \ref{exist}. This completes the proof of Theorem \ref{exist}.


{}


\begin{thebibliography}{References}

\bibitem{Cal}
M. C. Calderer.
\newblock {\em On the mathematical modeling of textures in polymeric liquid crystals}.
\newblock  Nematics (Orsay, 1990), 25�C36, NATO Adv. Sci. Inst. Ser. C Math. Phys. Sci., 332, Kluwer Acad. Publ., Dordrecht, 1991.

\bibitem{CC}
M. C. Calderer,  and C. Liu.
\newblock {\em Liquid crystal flow: dynamic and static configurations}.
\newblock  SIAM J. Appl. Math. 60, no. 6, 1925�C1949, 2002 (electronic).

\bibitem{CC1}
M. C. Calderer,  and C. Liu.
\newblock {\em Mathematical developments in the study of smectic A liquid crystals}.
\newblock  The Eringen Symposium dedicated to Pierre-Gilles de Gennes (Pullman, WA, 1998). Internat. J. Engrg. Sci. 38, no. 9-10, 1113�C1128, 2000.

\bibitem{CDLL}
M. C. Calderer, D.  Golovaty, F-H. Lin and C. Liu.
\newblock {\em Time evolution of nematic liquid crystals with variable degree of orientation}.
\newblock SIAM J. Math. Anal. 33, no. 5, 1033�C1047, 2002 (electronic).

\bibitem{CM} F. Crispo and P. Maremonti.
\newblock {\em An Interpolation Inequality in Exterior Domains}.
\newblock Rend. Sem. Mat. Univ. Padova, 112, 2004.

\bibitem{DQS}
M. Dai, J. Qing, and M. E. Schonbek.
\newblock {\em Regularity of Solutions to the Liquid Crystals Systems in $\mathbb{R}^2$ and $\mathbb{R}^3$}.
\newblock submitted to Nonlinear Analysis, preprint, 2011.

\bibitem{Er0} J. L. Ericksen.
\newblock {\em Conservation Laws for Liquid Crystals}.
\newblock Trans. Soc. Rheol. 5 (1961) 22 - 34.

\bibitem{Er1} J. L. Ericksen.
\newblock {\em Continuum Theory of Nematic Liquid Crystals}.
\newblock Res Mechanica 21 (1987) 381- 392.

\bibitem{EK} J. L. Ericksen, and D. Kinderlehrer, eds..
\newblock {\em Theory and Applications of Liquid Crystals}.
\newblock IMA Vol. 5, Springer-Verlag, New York, 1986.

\bibitem{Ev} L. C. Evans.
\newblock {\em  Partial Differential Equations}.
\newblock Graduate Studies in Mathematics, Vol. 19.

\bibitem{F} A. Friedman
\newblock {\em  Partial Differential Equations}.

\bibitem{JT} F.  Jiang, and Zhong Tan.
\newblock {\em Global Weak Solution to the Flow of Liquid Crystals System}.
\newblock Math. Meth. Appl. Sci. (32)2009, 2243-2266.

\bibitem{HKL} D. Kinderlehrer, F-H. Lin, and R. Hardt.
\newblock {\em Existence and partial regularity of static liquid crystal configurations}.
\newblock Comm. Math. Phys. 105, no. 4, 547�C570, 1986.

\bibitem{Kin} D. Kinderlehrer.
\newblock {\em Recent Developments in Liquid Crystal Theory}.
\newblock Frontiers in pure and applied mathematics, 151�C178, North-Holland, Amsterdam, 1991.

\bibitem{LS} O. A. Ladyzhenskaya, and V. A. Solonnikov.
\newblock {\em Linear and Quasilinear Equations of Parabolic Type}.
\newblock Transl. Math. Monographs, Vol. 23, AMS 1986.



\bibitem{Le0} F. M. Leslie.
\newblock {\em Some Contitutive Equations for liquid crystals}.
\newblock Arch Rational Mech Anal. 28 (1968) 265 - 283.

\bibitem{Le1} F. M. Leslie.
\newblock {\em Theory of flow phenomena in liquid crystals}.
\newblock Advances in Liquid Crystals, Vol 4 G. Brown ed., Academic Press, New York,
1979 1- 81.

\bibitem{LL1} F. Lin, and C. Liu.
\newblock {\em Existence of Solutions for the Ericksen-Leslie System}.
\newblock Arch. Rational Mech. Anal. 154(2000), 135-156.

\bibitem{LL} F. Lin, and C. Liu.
\newblock {\em Nonparabolic Dissipative Systems Modeling the Flow of Liquid Crystals}.
\newblock Communications on Pure and Applied Mathematics, Vol. XLVIII(1995), 501-537.

\bibitem{LL2} F. Lin, and C. Liu.
\newblock {\em Partial regularity of the dynamic system modeling the flow of liquid crystals. }.
\newblock Discrete Contin. Dynam. Systems 2, no. 1, 1�C22, 1996.

\bibitem{Lc} C. Liu.
\newblock {\em An Introduction to Mathematical Theories of Elastic Complex Fluids}.
\newblock Notes,2006.

\bibitem{Liu} X. Liu, and Z. Zhang.
\newblock {\em Existence of the Flow of Liquid Crystals System}.
\newblock Chinese Annals of Math. Series A, 30(1), 2009.




\bibitem{Sch} M. Schonbek.
\newblock {\em $L^2$ Decay for Weak Solutions of the Navier-Stokes Equations}.
\newblock Archive for Rational Mechanics and Analysis, Vol. 88, No. 3, 209-222, 1985.

\bibitem{Sch1} M. Schonbek.
\newblock {\em Large Time Behavior of Solutions to the Navier-Stokes Equations}.
\newblock Comm. in Partial Differential Equations, 11(7), 733-763, 1986.

\bibitem{Sch95}
M. E. Schonbek.
\newblock {\em Large Time Behavior of Solutions to Navier-Stokes Equations in $H^m$ Spaces}.
\newblock Comm. in P.D.E, 20(1995), No. 1 and 2, 103-117.

\bibitem{Sch94}
M. E. Schonbek and M. Wiegner.
\newblock {\em On the Decay of Higher-Order Norms of the Solutions of Navier-Stokes Equations}.
\newblock Proc. Royal Society of Edinburgh Sect. A 126 (1996), no.3, 677-685.

\bibitem{Sch2} M. Schonbek.
\newblock {\em Uniform Decay Rates for Parabolic Conservation Laws}.
\newblock Journal of Nonlinear Analysis, Vol. 10, No. 9, 943-956, 1986.


\bibitem{Wu} H. Wu.
\newblock {\em Long-time Behavior for Nonlinear Hydrodynamic System Modeling the Nematic Liquid Crystal Flows}.
\newblock Discrete Contin. Dyn. Syst., 26, no. 1, 379-396, 2010.


\end{thebibliography}
\end{document}